\documentclass[12pt,a4paper]{amsart}
\usepackage{amssymb}
\usepackage{mathrsfs}

\headsep=1cm \numberwithin{equation}{section}
\newtheorem{thm}{Theorem}[section]
\newtheorem{cor}[thm]{Corollary}
\newtheorem{lem}[thm]{Lemma}
\newtheorem{prop}[thm]{Proposition}
\theoremstyle{definition}
\newtheorem{defn}[thm]{Definition}
\theoremstyle{remark}

\numberwithin{equation}{section}
\newcommand{\To}{\longrightarrow}

\newcommand{\FSF}{\rm FSF }
\newcommand{\nat}{\mathbb N}
\newcommand{\m}{\frak m }
\newcommand{\p}{\frak p }

\newcommand\Supp{\operatorname{Supp}}
\newcommand\Ass{\operatorname{Ass}}

\newcommand\Spec{\operatorname{Spec}}

\newcommand\Hom{\operatorname{Hom}}
\newcommand\Ext{\operatorname{Ext}}
\newcommand\grade{\operatorname{grade}}

\newcommand\Max{\operatorname{Max}}
\begin{document}
\title[Faltings' local-global principle of local cohomology]{Faltings' local-global principle for the minimaxness of local cohomology modules}
\author{Mohammad Reza  Doustimehr and Reza Naghipour$^*$}
\address{Department of Mathematics, University of Tabriz, Tabriz, Iran;
and School of Mathematics, Institute for Research in Fundamental
Sciences (IPM), P.O. Box 19395-5746, Tehran, Iran.}
\email{m\b{ }doustimehr@tabrizu.ac.ir}
\email{naghipour@ipm.ir} \email {naghipour@tabrizu.ac.ir}

\thanks{ 
This research  was been in part supported by a grant from IPM.\\
$^*$Corresponding author: e-mail: {\it naghipour@ipm.ir} (Reza Naghipour)}%
%\subjclass{}%

%\date{}%
%\dedicatory{}%
%\commby{}%
% ----------------------------------------------------------------
\begin{abstract}
The concept of Faltings' local-global principle for the minimaxness of local cohomology modules
over a commutative Noetherian ring $R$ is introduced, and it is shown that this principle holds at level 2. We also
establish the same principle at all levels over an arbitrary commutative Noetherian ring of dimension not exceeding 3.
These generalize  the main results of Brodmann et al. in \cite{BRS}. Moreover, it is shown that if $M$ is a
finitely generated $R$-module, $\frak a$ an ideal of $R$  and $r$ a non-negative integer such that $\frak a^tH^i_{\frak a}(M)$
is skinny for all $i<r$ and for some positive integer $t$, then for any minimax submodule $N$ of  $H^r_{\frak a}(M)$,
the $R$-module $\Hom_R(R/\frak a,  H^r_{\frak a}(M)/N)$ is finitely generated. As a consequence, it follows that the associated
primes of $H^r_{\frak a}(M)/N$ are finite.  This generalizes the main results of Brodmann-Lashgari \cite{BL} and Quy \cite{Qu}.
\end{abstract}
\maketitle
% ----------------------------------------------------------------
\section{Introduction}
Throughout this paper, let $R$ denote a commutative Noetherian ring
(with identity) and $\frak a$ an ideal of $R$. For an $R$-module $M$, the
$i^{\rm th}$ local cohomology module of $M$ with support in $V(\frak a)$
is defined as:
$$H^i_{\frak a}(M) = \underset{n\geq1} {\varinjlim}\,\, {\rm Ext}^i_R(R/\frak a^n, M).$$
Local cohomology was first defined and studied by Grothendieck. We
refer the reader to \cite{BS} or \cite{Gr1} for more details about
local cohomology.  An important theorem in local cohomology is
Faltings' local-global principle for the finiteness dimension of
local cohomology modules \cite[Satz 1]{Fa1}, which states that for a
positive integer $r$,   the $R_{\frak p}$-module $H^i_{\frak
aR_{\frak p}}(M_{\frak p})$ is finitely generated for all $i\leq r$
and for all ${\frak p}\in  {\Spec}\,R$ if and only if
the $R$-module $H^i_{\frak a}(M)$ is finitely generated for all $i\leq r$.\\
Another formulation of Faltings' local-global principle,
particularly relevant for this paper, is in terms of the
generalization of the finiteness dimension $f_{\frak a}(M)$ of $M$
relative to $\frak a$, where
$$f_{\frak a}(M):=\inf\{i\in \Bbb{N}_0\,\,|\,\,H^i_{\frak a}(M)\,\,{\rm is}\,\,{\rm not}\,\,{\rm finitely}\,\,{\rm generated}
\}.\,\,\,\,\,\,\,\,\,\,\,\,\,\,\,\,\,\,\,\,\,\,\,\,\,(\dag)$$
 With the usual convention that the
infimum of the empty set of integers is interpreted as $\infty$.
For any non-negative integer $n$, the {\it $n^{\rm th}$ finiteness dimension} $f^n_{\frak a}(M)$ of $M$
relative to $\frak a$ is defined  by $$f^n_{\frak a}(M):=\inf\{f_{\frak aR_{\frak p}}(M_{\frak
p})\,\,|\,\,{\frak p}\in {\rm Supp}(M/\frak a M)\,\,{\rm and}\,\,  \dim R/{\frak p}\geq n\}.$$ Note that $f^n_{\frak a}(M)$ is either a positive
integer or $\infty$ and that $f^0_{\frak a}(M)=f_{\frak a}(M)$. The {\it $n^{\rm th}$ finiteness dimension} $f^n_{\frak a}(M)$ of $M$
relative to $\frak a$ has been introduced by Bahmanpour et al. in \cite{BNS1} and they showed that
$$f_{\frak a}^1(M)=\inf\{i\in\mathbb{N}: H_{\frak a}^i(M)\text{ is not minimax}\}.$$
We shall show that
 $$f_{\frak a}^1(M)=\inf\{i\in\mathbb{N}_0:{\frak a}^tH_{\frak a}^i(M)\text{ is not minimax for all }t\in\mathbb{N}\}.$$
This motivates to introduce the notion of the ${\frak b}$-minimaxness dimension $\mu_{\frak
a}^{\frak b}(M)$ of $M$ relative to ${\frak a}$, by
$$\mu_{\frak a}^{\frak b}(M)=\inf\{i\in\mathbb{N}:{\frak b}^tH_{\frak a}^i(M)\, \text {is not minimax for all
}t\in\mathbb{N}\},$$ where $\frak b$ is a second ideal of $R$. Note that $\mu_{\frak a}^{\frak
a}(M)=f_{\frak a}^1(M)$.

Recall that the $\frak b$-{\it finiteness dimension} of $M$ relative to $\frak a$ is defined by
\begin{eqnarray*}
f_{\frak a}^{\frak b}(M)&:=& \inf\{i\in \Bbb{N}_0\,\,|\,\,\frak b\not\subseteq {\rm Rad}(0:_{R}H^i_{\frak a}(M))\}\\
&=& \inf\{i\in \Bbb{N}_0\,\,|\,\,{\frak b}^nH^i_{\frak a}(M)\neq 0\,\,{\rm for}\,\,{\rm all}\,\,n\in \Bbb{N}\}.
\end{eqnarray*}

 Brodmann et al. in \cite{BRS} defined and studied the concept of the local-global principle for annihilation
 of local cohomology modules at level $r\in\mathbb{N}$ for the ideals $\frak a$ and $\frak b$ of $R$. We say that the local-global principle for the annihilation of local cohomology modules holds at level $r$ if for every choice of ideals ${\frak a}$, ${\frak b}$ of $R$ and every choice of finitely generated $R$-module $M$,
 it is the case that  $$f_{{\frak a}R_{\frak p}}^{{\frak b}R_{\frak p}}(M_{\frak p})>r \,\,\,\,\, \text{ for all } {\frak p}\in{\Spec}\, R \Longleftrightarrow f_{\frak a}^{\frak b}(M)>r.$$
It is shown in \cite{BRS} that the local-global principle for the annihilation of local cohomology modules holds at levels 1,2, over an arbitrary commutative Noetherian ring $R$ and at all  levels whenever $\dim R\leq 4$.

We say that the local-global principle for the minimaxness of local cohomology modules holds at level
$r\in\mathbb{N}$ if for every choice of ideals ${\frak a}$, ${\frak b}$ of
$R$ with ${\frak b}\subseteq {\frak a}$ and every choice of finitely generated $R$-module $M$, it is the case that
$$\mu_{{\frak a}R_{\frak p}}^{{\frak b}R_{\frak p}}(M_{\frak p})>r\,\,\,\,\, \text{ for all }
{\frak p}\in{\Spec}\, R \Longleftrightarrow \mu_{{\frak a}}^{{\frak
b}}(M)>r.$$

Our main result in Section 2 is to introduce the concept of Faltings' local-global principle for the minimaxness of local cohomology modules
over a commutative Noetherian ring $R$, and we show that this principle holds at level 2. We also
establish the same principle at all levels over an arbitrary commutative Noetherian ring of dimension not exceeding 3.
Our tools for proving the main result in Section 2 is the following:
\begin{thm}
Let $R$ be a  Noetherian ring  and let  ${\frak b}$ be a second ideal of $R$ such that ${\frak b}\subseteq{\frak a}$.  Let $M$ be a finitely generated $R$-module and let  $r$ be a positive integer such that  the local cohomology
modules $H_{\frak a}^0(M),\dots, H_{\frak a}^{r-1}(M)$ are ${\frak a}$-cofinite. Then
$$\mu_{{\frak a}R_{\p}}^{{\frak b}R_{\p}}(M_{\p})>r \text{ for all } {\p}\in \Spec R \Longleftrightarrow \mu_{{\frak a}}^{{\frak b}}(M)>r.$$
\end{thm}

Pursuing this point of view further we establish the following consequence of Theorem 1.1  which is an extension
of the results of Brodmann et al. in \cite[Corollary 2.3]{BRS} and Raghavan in \cite{Ra} for an arbitrary Noetherian ring.
\begin{cor}
Let $R$ be a  Noetherian ring, $M$ a finitely generated $R$-module and $\frak a, \frak b$ two ideals of $R$ such that ${\frak b}\subseteq{\frak a}$
and $\frak aM\neq M$. Set  $r\in\{1, \grade_M\frak a, f_{\frak a}(M), f_{\frak a}^1(M), f_{\frak a}^2(M)\}$. Then
$$\mu_{{\frak a}R_{\p}}^{{\frak b}R_{\p}}(M_{\p})>r \,\,\,\,  \text{ for all } {\p}\in \Spec R \Longleftrightarrow \mu_{{\frak a}}^{{\frak b}}(M)>r.$$
\end{cor}

In Section 3, we explore an interrelation between this principle and the Faltings' local-global principle for the annihilation of local
cohomology modules, and show that the local-global principle for the annihilation of local cohomology modules holds at level $2$ over $R$ and at all levels whenever $\dim R\leq 3$. These generalize and reprove the main results of Brodmann et al. in \cite{BRS}.

An $R$-module $L$ is said to be a ${\FSF}$ module if there is a finitely generated submodule $N$ of $L$ such that support of quotient module $L/N$ is a finite set. The class of ${\FSF}$ modules was introduced by Quy \cite{Qu} and he has given some properties of these modules.
Another main result in Section 3 is following:
 \begin{prop}
Let $R$ be a  Noetherian ring and $M$ a finitely generated $R$-module. Let $\frak a$ be an ideal of $R$ and $r$ a positive integer such that the $R$-modules ${\frak a}^tH_{\frak a}^0(M)$,\dots, ${\frak a}^tH_{\frak a}^{r-1}(M)$ are ${\rm FSF}$
for some $t\in\mathbb{N}_0$. Then, for any minimax submodule $N$ of $H^r_{\frak a}(M)$, the $R$-module $\Hom_R(R/{\frak a},H_{\frak a}^r(M)/N)$ is finitely generated and  the $R$-modules
$H_{\frak a}^0(M)$,\dots, $H_{\frak a}^{r-1}(M)$ are ${\frak a}$-cofinite.
\end{prop}
We will call a module {\it skinny or weakly Laskerian module}, if each of its homomorphic images has only finitely many associated primes (cf. \cite{DM} and \cite{Ro}). By using Proposition 1.3 we obtain the following corollary which is a generalization of the main results of Brodmann-Lashgari \cite{BL} and Quy \cite{Qu}.
\begin{cor}
Let $R$ be a  Noetherian ring and $M$ a finitely generated $R$-module. Let $\frak a$ be an ideal of $R$ and  $r$  a positive integer such that the $R$-modules
${\frak a}^tH_{\frak a}^0(M)$,\dots, ${\frak a}^tH_{\frak a}^{r-1}(M)$ are skinny. Then, for any minimax submodule $N$ of $H^r_{\frak a}(M)$, the set $\Ass_R H^r_{\frak a}(M)/N$ is finite.
\end{cor}

Throughout this paper, $R$ will always be a commutative Noetherian
ring with non-zero identity and $\frak a$ will be an ideal of $R$. Recall
that an $R$-module $L$ is called $\frak a$-{\it cofinite} if $\Supp L\subseteq V(\frak a)$ and ${\rm Ext}^j_R(R/\frak a, L)$ is finitely
generated for all $j\geq 0$. The concept of $\frak a$-cofinite modules
were introduced by Hartshorne \cite{Ha}.  An $R$-module $L$ is said to be {\it minimax}, if there exists a
finitely generated submodule $N$ of $L$, such that $L/N$ is
Artinian. The class of minimax modules was introduced by H.
Z\"{o}schinger \cite{Zo1} and he has given in \cite{Zo1, Zo2} many
equivalent conditions for a module to be minimax.
We shall use $\Max R$ to denote the set of all maximal
ideals of $R$. Also, for any ideal $\frak a$ of $R$, we denote
$\{\frak p \in {\rm Spec}\,R:\, \frak p\supseteq \frak a \}$ by
$V(\frak a)$. Finally, for any ideal $\frak{b}$ of $R$, the {\it
radical} of $\frak{b}$, denoted by ${\rm Rad}(\frak{b})$, is defined to
be the set $\{x\in R \,: \, x^n \in \frak{b}$ for some $n \in
\mathbb{N}\}$. For any unexplained notation and terminology we refer
the reader to \cite{BS} and \cite{Mat}.

\section{Local-global principle for minimaxness of local cohomology}
In this section we introduce the concept of Faltings' local-global principle for the minimaxness of local cohomology modules
over a commutative Noetherian ring $R$, and we show that this principle holds at level 2. We also
establish the same principle at all levels over an arbitrary commutative Noetherian ring of dimension not exceeding 3.
We begin with the following lemmas which are needed in this section.
%%%%%%%%%%%%%%%%%%%%%%%%%%%%%%%proposition0%%%%%%%%%%%%%%%%%%%%%%%%%%%%%%%%%%%%%%%%%%%%%%%%%%%%%%%%%%%%%%%%%%%%%%%%%%%%%%%%%%%%%%%%%%%%%%%%%%%%%%%%%%%

\begin{lem}\label{lem1}
Let $R$ be a  Noetherian ring, ${\frak a}$ an ideal of $R$, and $M$ an arbitrary $R$-module. Then ${\frak
a}M$ is minimax if and only if $M/(0:_M{\frak a})$ is minimax.
\end{lem}

\proof This follows easily from the definition.  \qed \\

The next lemma, which is a generalization of \cite[Lemma 9.1.2]{BS}, states that the $R$-modules $H_{\frak a}^0(M)$,\dots, $H_{\frak a}^{s-1}(M)$ are minimax if and only if there is an integer $t\in\mathbb{N}$ such that ${\frak a}^tH_{\frak a}^0(M)$,\dots, ${\frak a}^tH_{\frak a}^{s-1}(M)$ are minimax.
%%%%%%%%%%%%%%%%%%%%%%%%%%%%%%%%%%%%%%%%%%%%lemma%%%%%%%%%%%%%%%%%%%%%%%%%%%%%%%%%%%%%%%%%%%%%%%%%%%%%%%%%%%%%%%%%%%%%%%%%%%%%%%%%%%%%%%%%%%%%%%%%%%%%%
\begin{lem}\label{lem2}
Let $R$ be a  Noetherian ring, ${\frak a}$ an ideal of $R$, and $M$ a finitely generated $R$-module. Let $s$ be a positive integer. Then the following statements are equivalent:

{\rm(i)}  $H_{\frak a}^i(M)$ is minimax for all $i<s$;

{\rm (ii)} There exists  a positive integer $t$ such that ${\frak a}^tH_{\frak a}^i(M)$ is minimax for all $i<s$.
\end{lem}

\proof The implication ${\rm(i)}\Longrightarrow {\rm(ii)}$ is obviously true. In order to show ${\rm(ii)}\Longrightarrow {\rm(i)}$,
we proceed by induction on $s$. If $s=1$ there is nothing to
show. Suppose that $s>1$ and the case $s-1$ is settled. By inductive
hypothesis the $R$-module $H_{\frak a}^i(M)$ is minimax for all
$i<s-1$, and so it is enough for us to show that the $R$-module $H_{\frak
a}^{s-1}(M)$ is minimax. To this end, since by virtue of Lemma 2.1, the $R$-module $H_{\frak
a}^{s-1}(M)/(0:_{H_{\frak a}^{s-1}(M)}{\frak a}^t)$ is minimax,
there exists a finitely generated submodule $N$ of $H_{\frak
a}^{s-1}(M)$ such that the $R$-module $$H_{\frak a}^{s-1}(M)/N+(0:_{H_{\frak
a}^{s-1}(M)}{\frak a}^t),$$ is Artinian.  On the other hand in view of \cite[Theorem
2.3]{BN},  the $R$-module $(0:_{H_{\frak a}^{s-1}(M)}{\frak a}^t)$ is finitely generated, and hence
$N+ (0:_{H_{\frak a}^{s-1}(M)}{\frak a}^t)$ is also  finitely generated. Therefore  $H_{\frak a}^{s-1}(M)$ is
minimax, as required.\qed \\
%%%%%%%%%%%%%%%%%%%%%%%%%%%%%%%%%%%%%%%%%%%%corollary3%%%%%%%%%%%%%%%%%%%%%%%%%%%%%%%%%%%%%%%%%%%%%%%%%%%%%%%%%%%%%%%%%%%%%%%%%%%%%
\begin{cor}\label{cor3}
Let $R$ be a  Noetherian ring, ${\frak a}$ an ideal of $R$, and $M$ a finitely generated $R$-module.
 Let $f_{\frak a}^1(M)$ denote the 1-th finiteness dimension of $M$ relative to ${\frak a}$. Then
$$f_{\frak a}^1(M)=\inf\{i\in\mathbb{N}_0:{\frak a}^tH_{\frak a}^i(M)\text{ is not minimax for all
}t\in\mathbb{N}\}.$$
\end{cor}

\proof  The result follows immediately from \cite[Corollary 2.4]{BNS1} and Lemma \ref{lem2}. \qed \\

Now, we introduce the notion of the ${\frak b}$-{\it minimaxness dimension} $\mu_{\frak a}^{\frak b}(M)$ of $M$ relative to ${\frak a}$, as a
generalization of the $\frak b$-{\it finiteness dimension} $f_{\frak a}^{\frak b}(M)$ of $M$ relative to $\frak a$.
%%%%%%%%%%%%%%%%%%%%%%%%%%%%%%%definition4%%%%%%%%%%%%%%%%%%%%%%%%%%%%%%%%%%%%%%%%%%%%%%%%%%%%%%%%%%%%%%%%%%%%%%%%%%%%%%%%%%%%%%%%%%%%%%%%%%%%%
\begin{defn}
Let $M$ be a finitely generated module over a Noetherian ring $R$ and let ${\frak b}, \frak a$ be two ideals of $R$ such that ${\frak b}\subseteq{\frak a}$. We define the ${\frak b}$-{\it minimaxness
dimension} $\mu_{\frak a}^{\frak b}(M)$ of $M$ relative to ${\frak a}$ by $$\mu_{\frak a}^{\frak b}(M):=\inf\{i\in\mathbb{N}:{\frak
b}^tH_{\frak a}^i(M)\text{ is not minimax for all}\,\, t \in\mathbb{N}\}.$$  Note that, since $\Gamma_{\frak a}(M)$ is
minimax, we can write $$\mu_{\frak a}^{\frak b}(M):=\inf\{i\in\mathbb{N}_0:{\frak b}^tH_{\frak a}^i(M)\text{ is not
minimax for all }t\in\mathbb{N}\},$$ that $\mu_{\frak a}^{\frak b}(M)$ is
either a positive integer or $\infty$,
 and, by Corollary \ref{cor3},
 $\mu_{\frak a}^{\frak a}(M)=f_{\frak a}^1(M)$.
\end{defn}
We can also introduce the Faltings' local-global principle for the minimaxness of local cohomology modules which is a generalization of the
 Faltings' local-global principle for the annihilation of local cohomology modules.

%%%%%%%%%%%%%%%%%%%%%%%%%%%%%%%%%%%%%%%%%%%%definition5%%%%%%%%%%%%%%%%%%%%%%%%%%%%%%%%%%%%%%%%%%%%%%%%%%%%%%%%%%%%%%%%%%%%%%%%%%%%%
\begin{defn}
Let $R$ be a commutative Noetherian ring and let $r$ be a positive integer. We say that the  {\it Faltings' local-global principle for the
minimaxness of local cohomology modules} holds at level $r$ (over the
ring $R$) if,  for every choice of ideals ${\frak a}$, ${\frak b}$ of
$R$ with ${\frak b}\subseteq{\frak a}$ and for every choice of finitely generated $R$-module  $M$, it
is the case that $$\mu_{{\frak a}R_{\p}}^{{\frak b}R_{\p}}(M_{\p})>r
\text{ for all } {\p}\in \Spec R \Longleftrightarrow \mu_{{\frak a}}^{{\frak b}}(M)>r$$
\end{defn}
The following theorem plays a key role in the proof of the main result of this section.
%%%%%%%%%%%%%%%%%%%%%%%%%%%%%%%%%%%%%%%%%%%%%%%%%%%%%%theorem7%%%%%%%%%%%%%%%%%%%%%%%%%%%%%%%%%%%%%%%%%%%%%%%%%%%%%%%%%%%%%%%%%%%%%%%%%%%%%%%%%%%%
\begin{thm}\label{thm6}
Let $R$ be a  Noetherian ring  and let  $\frak a, {\frak b}$ be two ideals of $R$ such that ${\frak b}\subseteq{\frak a}$.  Let $M$ be a finitely generated $R$-module and let  $r$ be a positive integer such that  the local cohomology
modules $H_{\frak a}^0(M),\dots, H_{\frak a}^{r-1}(M)$ are ${\frak a}$-cofinite. Then
$$\mu_{{\frak a}R_{\p}}^{{\frak b}R_{\p}}(M_{\p})>r \text{ for all } {\p}\in \Spec R \Longleftrightarrow \mu_{{\frak a}}^{{\frak b}}(M)>r.$$
\end{thm}

\proof Let $i$ be an arbitrary  non-negative integer such that $i\leqslant r$. It
is sufficient for us to show that there is a non-negative integer $t_0$ such
that ${\frak b}^{t_0}H_{\frak a}^i(M)$ is minimax.  To do this, in view of  \cite[Theorem
2.3]{BN},  the set $\Ass_R({\frak b}^tH_{\frak a}^i(M))$ is
finite,  for all  $t\in\mathbb{N}_0$. Thus for all $t\in\mathbb{N}_0$, the set $\Supp {\frak
b}^tH_{\frak a}^i(M)$ is a closed subset of $\Spec R$ (in the
Zariski topology), and so the descending chain
$$\dots \supseteq \Supp( {\frak b}^tH_{\frak a}^i(M))\supseteq \Supp ({\frak b}^{t+1}H_{\frak a}^i(M))\supseteq \dots$$
is eventually stationary. Therefore there is a
non-negative integer $t_0$ such that for each $t\geqslant t_0$, $$\Supp ({\frak b}^tH_{\frak a}^i(M))=\Supp ({\frak b}^{t_0}H_{\frak a}^i(M)).$$
Let ${\m}$ be a maximal ideal of $R$. Since $\mu_{{\frak
a}R_{\p}}^{{\frak b}R_{\p}}(M_{\p})>r$, for all  ${\p}\in \Spec R$,  it follows that  there is an integer
$u\geqslant t_0$ such that $({\frak b}R_{\m})^uH_{{\frak
a}R_{\m}}^i(M_{\m})$ is minimax. Now, let ${\p}\in \Spec R$ be such
that ${\p}\subsetneqq{\m}$. Then it follows from
$$(({\frak b}R_{\m})^uH_{{\frak a}R_{\m}}^i(M_{\m}))_{{\frak p}R_{\m}}\cong ({\frak b}^uH_{\frak a}^i(M))_{\p}$$
and the definition of minimax modules that  $({\frak b}^uH_{\frak a}^i(M))_{\p}$  is a finitely generated  $R_{\p}$-module (note that
a module $L$ which is minimax has the property that the localization $L_{\frak p}$ is a finitely generated $R_{\frak p}$-module for each non-maximal prime ideal $\frak p$).
Now, as
 $({\frak b}^uH_{\frak a}^i(M))_{\p}$ is ${\frak a}R_{\p}$-torsion, there is
an integer $v\geq1$ such that $({\frak b}^{u+v}H_{\frak a}^i(M))_{\p}=0$, and
so  ${\p}\not\in \Supp {\frak b}^{t_0}H_{\frak a}^i(M)$. Therefore
$\Supp {\frak b}^{t_0}H_{\frak a}^i(M)\subseteq \Max R$.  Furthermore,  in view of hypothesis and  \cite[Theorem 2.3]{BN},
the $R$-module $\Hom_R(R/\frak a, {\frak b}^{t_0}H_{\frak a}^i(M))$ is finitely generated.  Thus, as
$\Hom_R(R/\frak a, {\frak b}^{t_0}H_{\frak a}^i(M))\subseteq \Max R$  it follows that  $\Hom_R(R/\frak a, {\frak b}^{t_0}H_{\frak a}^i(M))$
is Artinian. As ${\frak b}^{t_0}H_{\frak a}^i(M)$ is $\frak a$ -torsion, it yields from Melkersson's theorem \cite [Theorem 1.3]{Me1} that
${\frak b}^{t_0}H_{\frak a}^i(M)$ is  Artinian. Hence ${\frak b}^{t_0}H_{\frak a}^i(M)$ is minimax, as required.  \qed \\

%%%%%%%%%%%%%%%%%%%%%%%%%%%%%%%%%%%%%%%%%%%%corollary8%%%%%%%%%%%%%%%%%%%%%%%%%%%%%%%%%%%%%%%%%%%%%%%%%%%%%%%%%%%%%%%%%%%%%%%%%%%%%
\begin{cor}\label{cor7}
Let $R$ be a  Noetherian ring, $M$ a finitely generated $R$-module and  ${\frak a}$  an ideal of $R$ with  $\dim M/\frak aM \leq1$.  Let ${\frak b}$ be a second ideal of $R$ such that ${\frak b}\subseteq{\frak a}$. Then, for any integer $r$,
$$\mu_{{\frak a}R_{\p}}^{{\frak b}R_{\p}}(M_{\p})>r \text{ for all } {\p}\in \Spec R \Longleftrightarrow \mu_{{\frak a}}^{{\frak b}}(M)>r.$$
\end{cor}

\proof The assertion follows from \cite[Corollary 3.5]{BNS1} and Theorem 2.6.\qed \\

%%%%%%%%%%%%%%%%%%%%%%%%%%%%%%%%%%%%%%%%%%%%corollary8.1%%%%%%%%%%%%%%%%%%%%%%%%%%%%%%%%%%%%%%%%%%%%%%%%%%%%%%%%%%%%%%%%%%%%%%%%%%%%%
\begin{cor}\label{cor8}
The local-global principle (for the minimaxness of local cohomology modules) holds  over any (commutative Noetherian) ring $R$ with
$\dim R\leq 2$ at all levels $r\in \mathbb{N}$.
\end{cor}

\proof The result follows easily from  \cite[Corollary 5.2]{CGH} and Theorem 2.6.\qed \\

Our next corollary is a generalization of a result of Raghavan \cite{Ra}.

%%%%%%%%%%%%%%%%%%%%%%%%%%%%%%%%%%%%%%%%%%%%%%%%%%%%%%theorem9%%%%%%%%%%%%%%%%%%%%%%%%%%%%%%%%%%%%%%%%%%%%%%%%%%%%%%%%%%%%%%%%%%%%%%%%%%%%%%%%%%%%
\begin{cor}\label{cor9}
The local-global principle (for the minimaxness of local cohomology modules) holds at level $1$ (over any commutative Noetherian ring).
\end{cor}

\proof The assertion follows from  Theorem \ref{thm6}.\qed \\

Our next corollary is  a generalization  of a result of Brodmann et al. \cite{BRS}.

\begin{cor}\label{cor10}
Let $R$ be a  Noetherian ring and let  $\frak a,  \frak b$ be two ideals of $R$ such that ${\frak b}\subseteq{\frak a}$.  Let $M$ be a finitely generated $R$-module such  that $\frak a M\neq M$. Then
$$\mu_{{\frak a}R_{\p}}^{{\frak b}R_{\p}}(M_{\p})>\grade_M \frak a\,\,\,\,  \text{ for all } {\p}\in \Spec R \Longleftrightarrow \mu_{{\frak a}}^{{\frak b}}(M)>\grade_M \frak a.$$
\end{cor}

\proof The assertion follows from the definition of $\grade_M \frak a$ and Theorem \ref{thm6}.\qed \\

The next result is a generalization of Corollary 2.10.

\begin{cor}\label{cor11}
Let $R$ be a  Noetherian ring and let  $\frak a,  \frak b$ be two ideals of $R$ such that ${\frak b}\subseteq{\frak a}$.  Let $M$ be a finitely generated $R$-module, and that  $r\in\{f_{\frak a}(M), f_{\frak a}^1(M), f_{\frak a}^2(M)\}$. Then
$$\mu_{{\frak a}R_{\p}}^{{\frak b}R_{\p}}(M_{\p})>r \,\,\,\,  \text{ for all } {\p}\in \Spec R \Longleftrightarrow \mu_{{\frak a}}^{{\frak b}}(M)>r.$$
\end{cor}

\proof The assertion follows from  \cite[Theorems  2.3 and 3.2]{BNS1} and  Theorem \ref{thm6}.\qed \\ %%%%%%%%%%%%%%%%%%%%%%%%%%%%%%%%%%%%%%%%%%%%theorem9.1%%%%%%%%%%%%%%%%%%%%%%%%%%%%%%%%%%%%%%%%%%%%%%%%%%%%%%%%%%%%%%%%%%%%%%%%%%%%%

We are now ready to state and prove the  main theorem of this section, which shows that Faltings' local-global principle  for the minimaxness of local cohomology modules is valid at level $2$ over any commutative Noetherian $R$. This generalizes the main result of Brodmann et al. in \cite[Theorem 2.6]{BRS}.
\begin{thm}\label{thm9.1}
The local-global principle (for the minimaxness of local cohomology modules) holds over any (commutative Noetherian) ring $R$  at level $2$.
\end{thm}

\proof Let $M$ be a finitely generated $R$-module such that
$\mu_{{\frak a}R_{\p}}^{{\frak b}R_{\p}}(M_{\p})>2 \,\,\,\,  \text{
for all } {\p}\in \Spec R$. We must show that $\mu_{{\frak
a}}^{{\frak b}}(M)>2$. In view of Corollary 2.9, it is enough for us
to show that there exists a non-negative integer $t$ such that the
$R$-module ${\frak b}^tH_{\frak a}^2(M)$ is minimax. To do this, let
$M'=M/\Gamma_{\frak b}(M)$. We first show that the $R$-module $\Hom
_R(R/{\frak a}, H_{\frak a}^2(M'))$ is finitely generated. Then the
exact sequence
$$0\To \Gamma_{\frak b}(M) \To M \To M' \To 0,$$
induces the exact sequence
$$\hspace{12mm}H_{\frak a}^1(M)\To H_{\frak a}^1(M')\To H_{\frak a}^2(\Gamma_{\frak b}(M))\To H_{\frak a}^2(M)\To H_{\frak a}^2(M').\hspace{13mm}(\dag)$$
Next, since the set $\Ass_R H_{\frak a}^1(M))$ is finite, it follows from the proof of Theorem 2.6 that, for every $\p\in\Spec R$ with $\dim R/\p>0$, there exists a non-negative integer $u$ such that  $({\frak b}^{u}H_{\frak a}^1(M))_{\p}=0$. Furthermore, there  exists a non-negative integer $v$ such that
${\frak b}^vH_{\frak a}^i(\Gamma_{\frak b}(M))=0$ for all $i\geq0$.  Thus it follows from the exact sequence
obtained by the localization of  the exact sequence $(\dag)$  at the prime ideal ${\p}$  with $\dim R/\p>0$ and \cite[Lemma 9.1.1]{BS} that
 $({\frak b}R_{\p})^kH_{{\frak a}R_{\p}}^1(M'_{\p})=0$, for some integer $k\in\nat_0$.
Moreover, by  \cite[Lemma 2.1.1]{BS}, there exists $x\in{\frak b}$ which
is a non-zerodivisor on $M'$. Then $x^kH_{{\frak a}R_{\p}}^1(M'_{\p})=0$, and if we consider the long exact
sequence of local cohomology modules (with respect to ${\frak a}R_{\p}$) induced by the short exact sequence
$$0\To  M'_{\p}\stackrel{x^k}\To M'_{\p}\To M'_{\p}/ x^k
M'_{\p}\To 0,$$ we see that $H_{{\frak a}R_{\p}}^1(M'_{\p})$ is a
homomorphic image of $H_{{\frak a}R_{\p}}^0(M'_{\p}/x^k M'_{\p})$, and so it is a finitely generated $R_{\p}$-module, for all
$\p\in\Spec R$ with $\dim R/\p>0$.  It therefore follows from
\cite[Proposition 2.2]{BNS1} that $H_{\frak a}^1(M')$ is minimax; and
hence, by \cite[Theorem 2.3]{BN}, the $R$-module  $\Hom_R(R/{\frak a}, H_{\frak a}^2(M'))$ is finitely generated.

Since $\mu_{\frak a R_{\p}}^{\frak b R_{\p}}(M_{\p})> 2$, analogous
to the  proof of Theorem \ref{thm6}, for each $\p\in \Spec R$ with
$\dim R/\p>0$, there is $u_{\p}\in\nat_0$ such that $(\frak
b^{u_{\p}}H_{\frak a}^2(M))_{\p}=0$. Thus it follows from the exact
sequence obtained by the localization of  the exact sequence
$(\dag)$  at the prime ideal ${\p}$  with $\dim R/\p>0$ and
\cite[Lemma 9.1.1]{BS} that
 $(\frak b^{v_{\p}}H_{\frak a}^2(M'))_{\p}=0$, for some integer $v_{\p}\in\nat_0$.
 Since  $\Hom_R(R/\frak a, H_{\frak a}^2(M'))$ is finitely generated, it follows from the proof of Theorem 2.6  that there is an integer
$k\in\nat_0$ such that $\frak b^{k}H_{\frak a}^2(M')$ is Artinian,
and so $\Supp ({\frak b}^{k}H_{\frak a}^2(M'))$ is a finite subset
of ${\Max}R$. Now,  let $s:=v+k$. Then $\Supp ({\frak b}^sH_{\frak
a}^2(M))\subseteq \Supp({\frak b}^{k}H_{\frak a}^2(M'))\subseteq
\Max R$. Let $\Supp({\frak b}^sH_{\frak a}^2(M)):=\{\m_1,\dots,
\m_r\}$. Then by assumption, for each integer $j$ with $1\leqslant j
\leqslant r$, there is a non-negative integer $s_j\geqslant s$ such
that $({\frak b}^{s_j}H_{\frak a}^2(M))_{\m_j}$ is a minimax
$R_{\m_j}$-module. If we set $t=\max\{s_1,\dots,s_r\}$, then
$\Supp({\frak b}^tH_{\frak a}^2(M))\subseteq \{\m_1,\dots,\m_r\}$
and $({\frak b}^tH_{\frak a}^2(M))_{\m_j}$ is minimax, for all $j$
with $1\leqslant j \leqslant r$. Therefore, by \cite[Theorem
3.3]{Ba} and \cite[Proposition 2.2]{AM}, ${\frak b}^tH_{\frak
a}^2(M)$ is a minimax
$R$-module, as required. \qed \\

\begin{cor}
The local-global principle (for the minimaxness of local cohomology modules) holds over any (commutative Noetherian) ring $R$  with $\dim R\leq 3$.
\end{cor}
\proof  The assertion follows from Corollary 2.9, Theorem 2.12 and \cite[Exercise 7.1.7]{BS}. \qed \\

\section{Annihilation and associated primes of local cohomology modules}
The main goal of this section is to explore an interrelation between the Faltings' local-global principle for the minimaxness and annihilation of local cohomology modules, and to show the local-global principle for the annihilation of local cohomology modules holds at level $2$ over $R$ and at all levels whenever $\dim R\leq 3$. These  reprove  the main results of Brodmann et al. in \cite{BRS}. Moreover, it will be shown in this section that the subjects of the previous section can be used to prove a finiteness result about local cohomology module. In fact, we will generalize the main results of
Brodmann-Lashgari and Quy. The main result is Theorem 3.6.

The following  result describes a relation between the local-global principle for the minimaxness and  the annihilation (of
local cohomology modules) over a commutative Noetherian ring $R$.
\begin{prop}\label{prop100}
The local-global principle for the minimaxness (of local cohomology modules) implies the local-global principle for the annihilation (of
local cohomology modules) over any (commutative Noetherian) ring $R$.
\end{prop}

\proof  Let $r$ be a non-negative integer, and  suppose that the local-global principle for the
minimaxness of local cohomology modules holds at level $r$.  Let $M$ be a finitely generated $R$-module such that
$f_{\frak a R_{\p}}^{\frak b R_{\p}}(M_{\p})>r$ for all $\p\in\Spec R$. We
must show that $f_{\frak a}^{\frak b}(M)>r$. To this end, as $f_{\frak a R_{\p}}^{\frak b R_{\p}}(M_{\p})>r$,  it follows that $\mu_{\frak a R_{\p}}^{\frak b
R_{\p}}(M_{\p})>r$. Therefore, by hypothesis  $\mu_{\frak a}^{\frak b}(M)>r$; and  hence there
exists $t\in\nat$ such that ${\frak b}^tH_{\frak a}^i(M)$ is minimax for all
$i\leq r$. Thus  the set $\Ass_R({\frak b}^tH_{\frak a}^i(M))$ is finite.  Let
$\Ass_R ({\frak b}^tH_{\frak a}^i(M)):=\{\p_1,\dots,\p_s\}$. By assumption, for
each $1\leq j\leq s$, there is an integer $t_j\geq t$ such that
$({\frak b}^{t_j}H_{\frak a}^i(M))_{\p_j}=0$. Set $k:=\max\{t_1,\dots, t_s\}$.
Then $\Ass_R({\frak b}^kH_{\frak a}^i(M))=\emptyset$; and hence
${\frak b}^kH_{\frak a}^i(M)=0$. Therefore $f_{\frak a}^{\frak b}(M)>r$, as required. \qed \\

As a consequence of previous proposition and Theorem 2.12, the following corollary shows that the local-global principle for the
annihilation of local cohomology modules holds at level $2$ over $R$. This  reproves the main result of Brodmann et al. in \cite{BRS}.

\begin{cor}
The Local-global principle (for the annihilation of local cohomology modules) holds over any (commutative Noetherian) ring $R$  at level $2$.
\end{cor}
\proof The assertion follows from Proposition 3.1 and Theorem 2.12. \qed \\

In \cite{Qu}, P. H. Quy introduced the class of ${\FSF}$ modules and has given some properties of these modules. An $R$-module $L$ is said to be a ${\FSF}$ module if there is a finitely generated submodule $N$ of $L$ such that support of quotient module $L/N$ is a finite set. It is shown in \cite[Theorem 3.3]{Ba} that an  $R$-module $L$ is ${\FSF}$ if and only if it is skinny.
%%%%%%%%%%%%%%%%%%%%%%%%%%%%%%%%%%%%%%%%%%%%%%%%%%%%%%theorem10%%%%%%%%%%%%%%%%%%%%%%%%%%%%%%%%%%%%%%%%%%%%%%%%%%%%%%%%%%%%%%%%%%%%%%%%%%%%%%%%%%%%
\begin{prop}\label{prop10}
Let $R$ be a Noetherian ring,  $M$ a finitely generated $R$-module, $\frak a$ an ideal of $R$ and  $r$  a positive integer such that the $R$-modules ${\frak a}^tH_{\frak a}^0(M)$,\dots, ${\frak a}^tH_{\frak a}^{r-1}(M)$ are ${\rm FSF}$
for some $t\in\mathbb{N}_0$. Then the $R$-module $\Hom_R(R/{\frak a},H_{\frak a}^r(M))$ is finitely generated and  the $R$-modules
$H_{\frak a}^0(M)$,\dots, $H_{\frak a}^{r-1}(M)$ are ${\frak a}$-cofinite. In particular the set $\Ass_R H^r_{\frak a}(M)$ is finite.
\end{prop}

\proof Since  for each $i<r$, ${\frak a}^tH_{\frak a}^i(M)$ is ${\rm FSF}$, it follows that there is a finitely generated submodule
$N_i$ of ${\frak a}^tH_{\frak a}^i(M)$ such that  the set $\Supp({\frak a}^tH_{\frak a}^i(M)/N_i)$ is  finite, and so
 $\dim ({\frak a}^tH_{\frak a}^i(M)/N_i)\leq1$, for each $i<r$. Therefore,
for each ${\p}\in \Spec R$ with $\dim R/{\p}>1$, we have
$$({\frak a}R_{\p})^tH_{{\frak a}R_{\p}}^i(M_{\p})\cong({\frak a}^tH_{\frak a}^i(M))_{{\p}}\cong (N_i)_{\p}.$$
 Hence the $R_{\p}$-module $({\frak a}R_{\p})^tH_{{\frak a}R_{\p}}^i(M_{\p})$ is
finitely generated, for each $i<r$. Now, as $({\frak a}R_{\p})^tH_{{\frak a}R_{\p}}^i(M_{\frak p})$ is ${\frak a}R_{\p}$-torsion, so there
exists a non-negative integer $s$ such that $({\frak a}R_{\p})^{t+s}H_{{\frak a}R_{\p}}^i(M_{\p})=0$. Thus, by
\cite[Proposition 9.1.2]{BS}, $H_{{\frak a}R_{\p}}^i(M_{\p})$ is a
finitely generated  $R_{\p}$-module, for every ${\p}\in \Spec R$ with $\dim R/{\p}>1$ and for all $i<r$.
It therefore follows from \cite[Proposition 3.1]{BNS1} that  $\Hom_R(R/{\frak a}, H_{\frak a}^r(M))$ is finitely generated and the
$R$-modules $H_{\frak a}^0(M)$,\dots, $H_{\frak a}^{r-1}(M)$ are ${\frak a}$-cofinite.\qed \\
%%%%%%%%%%%%%%%%%%%%%%%%%%%%%%%%%%%%%%%%%%%%corollary11%%%%%%%%%%%%%%%%%%%%%%%%%%%%%%%%%%%%%%%%%%%%%%%%%%%%%%%%%%%%%%%%%%%%%%%%%%%%%
\begin{cor}\label{cor11}
Let $R$ be a Noetherian ring, $M$  a finitely generated $R$-module, $\frak a$ an ideal of $R$ and  $r$  a positive integer such that the $R$-modules
${\frak a}^tH_{\frak a}^0(M)$,\dots,${\frak a}^tH_{\frak a}^{r-1}(M)$ are skinny. Then  the $R$-modules $H_{\frak a}^0(M)$,\dots,
$H_{\frak a}^{r-1}(M)$ are ${\frak a}$-cofinite and for any  ideal $\frak b$ of $R$ with $\frak b\subseteq \frak a$,
$$\mu_{{\frak a}R_{\p}}^{{\frak b}R_{\p}}(M_{\p})>r\,\,\, \text{ for all } {\p}\in \Spec R \Longleftrightarrow \mu_{{\frak a}}^{{\frak b}}(M)>r.$$
\end{cor}

\proof Apply \cite[Theorem 3.3]{Ba}, Proposition \ref{prop10} and Theorem \ref{thm6}.\qed \\

%%%%%%%%%%%%%%%%%%%%%%%%%%%%%%%%%%%%%%%%%%%%corollary12%%%%%%%%%%%%%%%%%%%%%%%%%%%%%%%%%%%%%%%%%%%%%%%%%%%%%%%%%%%%%%%%%%%%%%%%%%%%%
\begin{cor}\label{cor12}
Let $R$ be a Noetherian ring, $M$  a finitely generated $R$-module, $\frak a$ an ideal of $R$ and  $r$  a positive integer such that the $R$-modules
${\frak a}^tH_{\frak a}^0(M)$,\dots, ${\frak a}^tH_{\frak a}^{r-1}(M)$ have finite support, for some $t\in\mathbb{N}_0$. Then
the $R$-modules $H_{\frak a}^0(M),\dots, H_{\frak a}^{r-1}(M)$ are ${\frak a}$-cofinite and for any  ideal $\frak b$ of $R$ with $\frak b\subseteq \frak a$,
$$\mu_{{\frak a}R_{\p}}^{{\frak b}R_{\p}}(M_{\p})>r\,\,\, \text{ for all } {\p}\in \Spec R \Longleftrightarrow \mu_{{\frak a}}^{{\frak b}}(M)>r.$$
\end{cor}

\proof Since the set $\Supp ({\frak a}^tH_{\frak a}^i(M))$ is finite, for every $i<r$, it follows that the $R$-module ${\frak a}^tH_{\frak a}^i(M)$
is ${\FSF}$ for each $i<r$.   Now the assertion follows from Corollary 3.4. \qed \\

The following theorem, which is one of our main results, generalizes the main results of Brodmann-Lashgari \cite{BL} and Quy \cite{Qu}.

\begin{thm}
Let $R$ be a Noetherian ring, $M$  a  finitely generated $R$-module
 and $\frak a$ an ideal of $R$. Let $r$ be a non-negative integer such that $\frak a^tH^{i}_{I}(M)$ is
${\FSF}$ for all $i<r$ and for some $t\in\mathbb{N}_0$.  Then, for any minimax submodule $N$ of
$H^{r}_{\frak a}(M)$,   the $R$-module $\Hom_{R}(R/\frak a, H^{r}_{\frak a}(M)/N)$ is finitely generated. In particular,
the set $\Ass_{R}(H^{r}_{\frak a}(M)/N)$ is finite.
\end{thm}
\proof In view of Proposition 3.3, $\Hom_{R}(R/\frak a, H^{r}_{\frak a}(M))$ is
finitely generated. On the other hand, according to  Melkersson
\cite[Proposition 4.3]{Me2}, $N$ is $\frak a$-cofinite. Now, the exact
sequence $$0 \longrightarrow N \longrightarrow H^{r}_{\frak a}(M)
\longrightarrow H^{r}_{\frak a}(M)/N \longrightarrow 0$$ induces the
following exact sequence,$$\Hom_{R}(R/\frak a,H^{r}_{\frak a}(M))
\longrightarrow \Hom_{R}(R/\frak a,H^{r}_{\frak a}(M)/N) \longrightarrow
\Ext^{1}_{R}(R/\frak a,N).$$ Consequently $\Hom_{R}(R/\frak a, H^{r}_{\frak a}(M)/N)$ is finitely generated, as
required.\qed\\

In  \cite[Theorem 2.3]{BN}, Bahmanpour and Naghipour showed that if  $H_{\frak a}^0(M)$,\dots, $H_{\frak a}^{r-1}(M)$ are
minimax, then $H_{\frak a}^0(M)$,\dots, $H_{\frak a}^{r-1}(M)$ are ${\frak a}$-cofinite. The following corollary provides a slight generalization of \cite[Theorem 2.3]{BN}.
%%%%%%%%%%%%%%%%%%%%%%%%%%%%%%%%%%%%%%%%%%%%corollary13%%%%%%%%%%%%%%%%%%%%%%%%%%%%%%%%%%%%%%%%%%%%%%%%%%%%%%%%%%%%%%%%%%%%%%%%%%%%%
\begin{cor}\label{cor13}
Let $R$ be a Noetherian ring, $M$  a finitely generated $R$-module, $\frak a$ an ideal of $R$ and  $r$  a positive integer such that the $R$-modules
${\frak a}^tH_{\frak a}^0(M)$,\dots, ${\frak a}^tH_{\frak a}^{r-1}(M)$ are skinny. Then the set $\Ass_R H^r_{\frak a}(M)$ is finite.
\end{cor}
\proof Apply \cite[Theorem 3.3]{Ba} and Theorem 3.6.\qed \\

\begin{center}
{\bf Acknowledgments}
\end{center} The authors would like to thank Professor Hossein Zakeri for his reading of the first draft and valuable discussions.
We also would like to thank from the Institute for Research in Fundamental Sciences (IPM), for its financial support.
% ----------------------------------------------------------------

\end{document}